\documentclass[12pt,leqno]{article}
\usepackage[english]{babel}

\usepackage{amsfonts}
\usepackage{amssymb}
\usepackage{eucal}
\usepackage{amsxtra}
\usepackage{color}
\addtocounter{section}{-1} \numberwithin{equation}{section}
\newtheorem{theorem}{Theorem}[section]

\newtheorem{lemma}[theorem]{Lemma}

\newcommand{\qed}{{$\hfill \Box$}}

\begin{document}

\title{REAL HYPERSURFACES OF NON - FLAT COMPLEX PLANES IN TERMS OF THE JACOBI OPERATORS}
\bigskip
\author{Th. Theofanidis,  Ph. J. Xenos}

\date{}

\maketitle
\begin{minipage}{400pt}
\begin{center}
Mathematics Division, School of Technology,\\
Aristotle University of Thessaloniki,\\
Thessaloniki, 54124, Greece.\\
email  :theotheo@gen.auth.gr, fxenos@gen.auth.gr.
\end{center}

\begin{abstract}\hspace{-18pt}The aim of the present paper is the study of
some classes of real hypersurfaces equipped with the condition $\phi
l = l \phi$, $l = R( . , \xi)\xi.$
\end{abstract}
\
\end{minipage}

\section{Introduction.}
An n - dimensional Kaehlerian manifold of constant holomorphic
sectional curvature c is called complex space form, which is denoted
by $M_{n}(c)$. The complete and simply connected complex space form
is a projective space $\mathbb{C}P^{n}$ if $c > 0$, a hyperbolic
space $\mathbb{C}H^{n}$ if $c < 0$, or a Euclidean space
$\mathbb{C}^{n}$ if $c = 0$. The induced almost contact metric
structure of a real hypersurface M of $M_{n}(c)$ will be denoted by
($\phi, \xi, \eta, g$).

Real hypersurfaces in $\mathbb{C}P^{n}$ which are homogeneous, were
classified by R. Takagi (\cite{Takagi}). J. Berndt (\cite{Berndt})
classified  real hypersurfaces with principal structure vector
fields in $\mathbb{C}H^{n}$, which are divided into the model spaces
$A_{0}$, $A_{1}$, $A_{2}$ and $B$. Another class of real
hypersurfaces were studied by Okumura \cite{Okumura}, and Montiel
and Romero \cite{Montiel Romero}, who proved respectively the
following theorems.
\begin{theorem}
Let M be a real hypersurface of  $\mathbb{C}P^{n}$, $n \geq 2$. If
it satisfies
$$ g((A \phi - \phi A)X, Y) = 0$$
for any vector fields X and Y, then M
is a tube of radius r over one of the following Kaehlerian submanifolds:\\
\
$(A_{1})$ a hyperplane $\mathbb{C}P^{n-}$, where $0 < r <\frac{\pi}{2}$,\\
$(A_{2})$ a totally geodesic $\mathbb{C}P^{k}$($0 < k \leq n -
2$),where $0 < r <\frac{\pi}{2}$.
\end{theorem}

\begin{theorem}
Let M be a real hypersurface of  $\mathbb{C}H^{n}$, $n \geq 2$. If
it satisfies
$$ g((A \phi - \phi A)X, Y) = 0$$
for any vector fields X and Y, then M is locally congruent to
one of the following:\\
$(A_{0})$ a self - tube, that is, horosphere,\\
$(A_{1})$ a geodesic hypershere or a tube over a hyperplane $\mathbb{C}H^{n - 1}$,\\
$(A_{2})$ a tube over a totally geodesic $\mathbb{C}H^{k}$ ($1 \leq
k \leq n - 2$).
\end{theorem}

Real hypersurfaces of type $A_{1}$ and $A_{2}$ in $\mathbb{C}P^{n}$
and of type $A_{0}$, $A_{1}$ and $A_{2}$ in $\mathbb{C}H^{n}$ are
said to be hypersurfaces of \emph{type  A} for simplicity.

A Jacobi field along geodesics of a given Riemannian manifold (M, g)
plays an important role in the study of differential geometry. It
satisfies a well known differential equation which inspires Jacobi
operators. For any vector field $X$, the Jacobi operator is defined
by $R_{X}$: $R_{X}(Y) = R(Y, X))X$, where $R$ denotes the curvature
tensor and $Y$ is a vector field on M. $R_{X}$ is a self - adjoint
endomorphism in the tangent space of M, and is related to the Jacobi
differential equation, which is given by
$\nabla_{\acute{\gamma}}(\nabla_{\acute{\gamma}}Y) + R(Y,
\acute{\gamma})\acute{\gamma} = 0$ along a geodesic $\gamma$ on $M$,
where $\acute{\gamma}$ denotes the velocity vector along $\gamma$ on
$M$.

In a real hypersurface $M$ of a complex space form $M_{n}(c)$, $c
\neq 0$,
 the Jacobi operator on $M$ with respect to the structure vector field $\xi$,
 is called the structure Jacobi operator and is denoted by
$R_{\xi}(X) = R(X, \xi,)\xi$.

Many authors have studied real hypersurfaces from many points of
view. Certain authors have studied real hypersurfaces under the
condition $\phi l = l \phi$, equipped with one or two additional
conditions (\cite{Cho Ki}, \cite{Ki}, \cite{Ki Kim Lee}, \cite{Ki
Lee Lee} \cite{Ki Nagai Takagi}). Any such hypersurface is locally
congruent to a model space of type $A$.

In 2006, Ortega, Perez and Santos \cite{Ortega Perez Santos}in
showed the nonexistence of real hypersurfaces in nonflat complex
space forms whose structure Jacobi operator is parallel. In 2007,
Ki, Perez, Santos and Suh in \cite{Ki Perez Santos} classified real
hypersurfaces in complex space forms with $\xi$- parallel Ricci
tensor and structure Jacobi operator. Recently J. T. Cho and U - H
Ki in \cite{Cho Ki1} classified the real hypersurfaces whose
structure Jacobi operator is symmetric along the Reeb flow $\xi$ and
commutes with the shape operator A.

 In the present paper we consider a weaker condition $\nabla_{\xi}l = \mu\xi$
 where $\mu$ is a function of
class $C^{1}$ on $M$ and classify these hypersurfaces
satisfying $\phi l = l\phi$. Namely we prove:
\begin{theorem}
Let M be a real hypersurface of a complex plane $M_{2}(c)$ $(c
\neq 0)$, satisfying
$\phi l = l \phi$. If $lA = Al$ on $ker(\eta)$ or on span$\{\xi\}$ then M is locally congruent to a model space of type A.\\
\end{theorem}

\begin{theorem} Let M be a real hypersurface of a complex plane
$M_{2}(c)$ $(c \neq 0)$, satisfying
$\phi l = l \phi$. If  $\nabla_{\xi}l = \mu \xi$ on $ker(\eta)$ or on span$\{\xi\}$, then M is locally congruent to a model space of type A.\\
\end{theorem}

\section{Preliminaries.}
 Let $M_{n}$ be a Kaehlerian manifold of real
dimension 2n, equipped with an almost complex structure J and a
Hermitian metric tensor G. Then for any vector fields X and Y on
$M_{n}(c)$, the following
relations hold:\\
\begin{center}
$J^{2}X = -X$, \hspace{20pt}$G(JX, JY) = G(X, Y)$,\hspace{20pt}
$\widetilde{\nabla}J = 0$\end{center} where $\widetilde{\nabla}$
denotes the Riemannian connection of G of
$M_{n}$.\\
\
\\
Now, let $M_{2n-1}$ be a real (2n-1)-dimensional hypersurface of
$M_{n}(c)$, and denote by N a unit normal vector field on a
neighborhood of a point in $M_{2n-1}$ (from now on we shall write
\emph{M} instead of $M_{2n-1}$). For any vector field  X tangent to
M we have $JX = \phi X + \eta(X)N$, where $\phi X$ is the tangent
component of $JX$, $\eta(X)N$ is the normal component, and
\begin{center}
$\xi = - JN$, \hspace{30pt} $\eta(X) = g(X, \xi)$, \hspace{30pt}$g =
G|_{M}$.
\end{center}

By properties of the almost complex structure J, and the definitions
of $\eta$ and g, the following relations hold (\cite{Blair}):

\begin{center}
\begin{flushleft}
(1.1)\hspace{70pt}$\phi^{2} = - I + \eta\otimes\xi, \hspace{20pt}
\eta\circ\phi = 0 ,
    \hspace{20pt} \phi\xi = 0, \hspace{10pt}
 \eta(\xi) = 1 $\newline
\ \\ (1.2)\hspace{20pt}$g(\phi X, \phi Y) = g(X, Y) -
\eta(X)\eta(Y), \hspace{20pt} g(X, \phi Y) = - g(\phi X,
Y)$.\end{flushleft}

\end{center}
The above relations define an \emph{almost contact metric structure}
on M which is denoted by $(\phi, \xi, g, \eta)$. When an almost
contact metric structure is defined on M, we can define a local
orthonormal basis $\{V_{1}, V_{2}, . . . V_{n - 1}, \phi V_{1}, \phi
V_{2}, . . . \phi V_{n - 1}, \xi\}$, called a $\phi - basis$.
Furthermore, let A be the shape operator in the direction of N, and
denote by $\nabla$ the Riemannian connection of g on M. Then, A is
symmetric and the following equations are
satisfied:\\
\begin{center}
\begin{flushleft}
(1.3) \hspace{40pt} $\nabla_{X}\xi = \phi AX  ,\hspace{40pt}
(\nabla_{X}\phi)Y = \eta(Y)AX -
    g(AX, Y)\xi$.\end{flushleft}
\end{center}

As the ambient space $M_{n}(c)$ is of constant holomorphic sectional
curvature c, the equations of Gauss and Godazzi are respectively
given by:

\begin{center}
(1.4) \hspace{20pt}$R(X, Y)Z = \frac{c}{4}[g(Y, Z) X - g(X, Z)Y +
g(\phi Y, Z)\phi X - g(\phi
X, Z)\phi Y$\\
\ \\
\hspace{20pt}$ - 2g(\phi X, Y)\phi Z] + g(AY, Z)AX - g(AX,
Z)AY$,\end{center} \
\begin{center}
(1.5)\hspace{60pt}$ (\nabla_{X}A)Y - (\nabla_{Y}A)X =
\frac{c}{4}[\eta(X)\phi Y - \eta(Y)\phi X - 2g(\phi X, Y)\xi]$.
\end{center}
The tangent space $T_{p}M$, for every point $p\in M$, is decomposed
as following:
\begin{center}
$T_{p}M = ker(\eta)^{\bot} \oplus ker(\eta)$
\end{center}
where $ker(\eta)^{\bot} = span\{\xi\}$ and $ker(\eta)$ is defined as following:\\
\begin{center}
$ker(\eta) = \large{\{}
X \in T_{p}M: \eta(X) = 0\large{\}}$ \\
\end{center}
\
\\
Based on the above decomposition, by virtue of (1.3), we decompose
the vector field $A\xi$ in the following way:
\begin{flushleft}
(1.6) \hspace{100pt}$A\xi = \alpha \xi + \beta U$\end{flushleft}

\begin{flushleft}
where $\beta = |\phi\nabla_{\xi}\xi|$ and $U = -\frac{1}{\beta}
\phi\nabla_{\xi}\xi \in ker(\eta)$, provided that $\beta \neq 0$.
\end{flushleft}
If the vector field $A\xi$ is expressed as $A\xi = \alpha \xi$, then
$\xi$ is called a principal vector field.

Finally differentiation will be denoted by ( ). All manifolds of
this paper are assumed to be connected and of class $C^{\infty}$.

\section{Hypersurfaces satisfying $\phi l = l \phi$}
In the study of real hypersurfaces of a complex space form
$M_{n}(c)$, $c \neq 0$, it is a crucial condition that the structure
vector field $\xi$ is principal.  The purpose of this paragraph is
to prove this condition.\\
Let V be the open subset of points p of M, where $\alpha\neq0$ in a
neighborhood of p and $V_{0}$ be the open subset of points p of M
such that $\alpha = 0$ in a neighborhood of p. Since $\alpha$ is a
smooth function on M, then $V\cup V_{0}$ is an open and dense subset
of M.

\begin{lemma}
Let M be a real hypersurface of a complex plane $M_{2}(c)$ $(c
\neq 0)$, satisfying $\phi l = l \phi$ on $ker(\eta).$ Then, $\beta
= 0$ on $V_{0}$.
\end{lemma}
Proof.\\
From (1.6) we have $A\xi = \beta U$ on $V_{0}$. Then (1.4) for $X =
U$ and $Y = Z= \xi$ yields
$$lU = \frac{c}{4}U + g(A\xi, \xi)AU - g(AU, \xi)A\xi = \frac{c}{4}U
- g(U, A\xi)A\xi = (\frac{c}{4} - \beta^{2})U \Rightarrow$$
$$\phi l U = (\frac{c}{4} - \beta^{2})\phi U. $$
In the same way, from (1.4) for $X = \phi U$, $Y = Z = \xi$ we
obtain $l \phi U = \frac{c}{4} \phi U.$
The last two equations yield $\beta = 0.$ \qed\\
\emph{\textbf{REMARK}}\\
We have proved that on $V_{0}$, $A\xi = 0 \xi$ i.e. $\xi$ is a
principal vector field on $V_{0}$. Now we define on V the set $V'$
of points p where $\beta \neq 0$ in a neighborhood of p and the set
$V''$ of points p where $\beta = 0$ in a neighborhood of p.
Obviously $\xi$ is principal on $V''$. In what follows we study the
open subset $V'$ of M.

\begin{lemma}
Let M be a real hypersurface of a complex plane $M_{2}(c)$
($c\neq 0$), satisfying $\phi l = l \phi$ on $ker(\eta)$ . Then the following relations hold on $V'$.\\
\end{lemma}
\begin{equation}
AU = \Big(\frac{\gamma}{\alpha}-\frac{c}{4\alpha} +
\frac{\beta^{2}}{\alpha}\Big)U + \beta \xi, \hspace{50pt}A \phi U =
(\frac{\gamma}{\alpha} - \frac{c}{4\alpha})\phi U
\end{equation}

\begin{equation}
\nabla_{\xi}\xi = \beta \phi U , \hspace{20pt}  \nabla_{U}\xi =
\Big(\frac{\gamma}{\alpha} - \frac{c}{4\alpha} +
\frac{\beta^{2}}{\alpha}\Big) \phi U ,\hspace{20pt} \nabla_{\phi
U}\xi = - (\frac{\gamma}{\alpha} - \frac{c}{4\alpha}) U
\end{equation}

\begin{equation}
\nabla_{\xi}U = \kappa_{1} \phi U , \hspace{20pt} \nabla_{U}U =
\kappa_{2} \phi U , \hspace{20pt} \nabla_{\phi U}U = \kappa_{3} \phi
U + (\frac{\gamma}{\alpha} - \frac{c}{4\alpha})\xi
\end{equation}

\begin{equation}
\nabla_{\xi}\phi U = -\kappa_{1} U - \beta \xi,
\hspace{10pt}\nabla_{U}\phi U = -\kappa_{2} U -
\Big(\frac{\gamma}{\alpha} - \frac{c}{4\alpha} +
\frac{\beta^{2}}{\alpha}\Big)\xi, \nabla_{\phi U}\phi U =
-\kappa_{3} U
\end{equation}
where $\kappa_{1}$, $\kappa_{2}$, $\kappa_{3}$ are smooth functions on V.\\
\
\\
Proof.\\
 Using (1.4) with $Y = Z = \xi$, and from the definition of $lX$, we get \
\\
\begin{equation}
    lX = \frac{c}{4}[X - \eta(X)\xi] + \alpha AX - g(AX, \xi)A\xi
\end{equation}
which, for X = U yields
\begin{equation}
    lU = \frac{c}{4}U + \alpha AU - \beta A\xi.
\end{equation}
The inner products of (2.6) with U and $\phi U$ yield respectively
\begin{equation}
    g(AU, U) = \frac{\gamma}{\alpha} - \frac{c}{4\alpha} + \frac{\beta^{2}}{\alpha},
\end{equation}
\begin{equation}
g(AU, \phi U) = \frac{1}{a}g(lU, \phi U).
\end{equation}
where $\gamma = g(lU, U).$\\
 The second relation
of (1.2) for $X = U, Y = lU$, the condition $\phi l =l \phi$ and the
symmetry of the operator $l$ imply: $g(lU, \phi U) = 0 $. The last
two equations imply
\begin{equation}
g(AU, \phi U) = 0.
\end{equation}
Taking the inner product of (2.6) and $\xi$ we get:
\begin{equation}
g(AU, \xi) = \beta.
\end{equation}
From relations (2.7), (2.9) and (2.10), we obtain the first of
(2.1).\\
Because of (1.6) and the symmetry of A we obtain
\begin{equation}
g(A \phi U, \xi) = 0.
\end{equation}
Putting $X = \phi U$ in (2.5)  we have :
\begin{equation}
l\phi U = \frac{c}{4}\phi U + \alpha A\phi U.
\end{equation}
Since $\gamma = g(l\phi U, \phi U)$, taking the inner product of
(2.12) with $\phi U$, we have
\begin{equation}
g(A \phi U,  \phi U) = \frac{\gamma}{\alpha} -\frac{c}{4\alpha}.
\end{equation}
Using (2.9), (2.11) and (2.13), we obtain the second of (2.1).\\
\ \\
If we put $X = \xi$,  $X = U$, $X = \phi U$ in  (1.3) respectively
and by making use of (1.6)
and (2.1), we obtain (2.2).\\
It is well known that:
\begin{equation}
X g(Y, Z) = g(\nabla_{X}Y, Z) +   g(Y, \nabla_{X}Z)
\end{equation}
\
\\
Now if we use (2.2) and (2.14), it is easy to verify that:
$g(\nabla_{\xi}U, U) = 0 = g(\nabla_{\xi}U, \xi) $. So if we put
$g(\nabla_{\xi}U, \phi U) = \kappa_{1}$, we have $\nabla_{\xi}U =
\kappa_{1}\phi U$.  From (2.2) and (2.14), we obtain $g(\nabla_{U}U,
U) = g(\nabla_{U}U, \xi) = 0$. Therefore, putting $g(\nabla_{U}U,
\phi U) = \kappa_{2}$, we have $\nabla_{U}U = \kappa_{2} \phi U$.
 By use of (2.14) and (2.2) we have that $g(\nabla_{\phi U}U, U ) = 0$
 and $g(\nabla_{\phi U}U, \xi ) = \frac{\gamma}{\alpha} - \frac{c}{4\alpha}$. Then if we set $g(\nabla_{\phi U}U, \phi U
 ) = \kappa_{3}$, we get
$\nabla_{\phi U}U = \kappa_{3} \phi U + (\frac{\gamma}{\alpha} -
\frac{c}{4\alpha})\xi$.
 In the
same way, using (2.14), (2.2) and ( 2.3) we obtain (2.4).\qed

\begin{lemma}
Let M be a real hypersurface of a complex plane $M_{2}(c)$ $(c
\neq 0)$, satisfying $\phi l = l \phi$ on $ker(\eta)$. Then on $V'$
we have $\gamma \neq \frac{c}{4}$.\end{lemma}
 Proof.\\
Putting $X = U$, $Y = \xi$ in (1.5), we obtain $(\nabla_{U}A)\xi -
(\nabla_{\xi}A)U = - \frac{c}{4}\phi U.$ Combining the last equation
with (1.6), (2.1), (2.2), (2.3) and (2.4) respectively, it follows :
\vspace{-15pt}
\begin{center}
$$\Big[(U\alpha) - (\xi\beta) \Big]\xi + \Big[(U\beta) - (\xi\big(\frac{\gamma}{\alpha}-\frac{c}{4\alpha} +
\frac{\beta^{2}}{\alpha}  \big))\Big]U +$$\\
$$\Big[\gamma - \frac{c}{4} + \kappa_{2}\beta -(\frac{\gamma}{\alpha} - \frac{c}{4\alpha})\big(\frac{\gamma}{\alpha} - \frac{c}{4\alpha}
+ \frac{\beta^{2}}{\alpha}  \big) -
\kappa_{1}\frac{\beta^{2}}{\alpha} \Big]\phi U = - \frac{c}{4}\phi
U.$$
\end{center}
The last equation because of the linear independency of $U$, $\phi U
$ and $\xi$, yields
\begin{equation}
\big(U\alpha\big) = \big(\xi\beta \big),
\end{equation}
\begin{equation}
\big(U\beta \big) = \Big(\xi\big(\frac{\gamma}{\alpha} -
\frac{c}{4\alpha} + \frac{\beta^{2}}{\alpha}) \Big),
\end{equation}
\begin{equation}
\gamma + \kappa_{2}\beta -(\frac{\gamma}{\alpha} -
\frac{c}{4\alpha})\big(\frac{\gamma}{\alpha} - \frac{c}{4\alpha} +
\frac{\beta^{2}}{\alpha}  \big) -
\kappa_{1}\frac{\beta^{2}}{\alpha}= 0.
\end{equation}
In the same way, putting $X = \phi U$, $Y = \xi$ in (1.5) we obtain
$(\nabla_{\phi U}A)\xi - (\nabla_{\xi}A)\phi U = \frac{c}{4}U.$
Combining the last equation with (1.6), (2.1), (2.2), (2.3) and
(2.4), we have
\begin{equation}
\big(\phi U \beta \big) +(\frac{\gamma}{\alpha}-
\frac{c}{4\alpha})\Big(\frac{\gamma}{\alpha} - \frac{c}{4\alpha} +
\frac{\beta^{2}}{\alpha}\Big) - \kappa_{1}\frac{\beta^{2}}{\alpha} -
\beta^{2} - \gamma= 0,
\end{equation}
\begin{equation}
\kappa_{3}\beta = \xi\big(\frac{\gamma}{\alpha} -
\frac{c}{4\alpha}\big),
\end{equation}
\begin{equation}
\big(\phi U \alpha \big) + 3\beta(\frac{\gamma}{\alpha} -
\frac{c}{4\alpha}) - \kappa_{1}\beta - \alpha\beta = 0.
\end{equation}
Similarly, putting $X = U$, $Y = \phi U$ in (1.5), we get
$(\nabla_{U}A)\phi U - (\nabla_{\phi U}A)U = - \frac{c}{2}\xi ,$
which, by use of (1.6), (2.1), (2.2), (2.3) and (2.4), implies :
\begin{equation}
-\kappa_{2}\frac{\beta^{2}}{\alpha} - 3\beta(\frac{\gamma}{\alpha}
-\frac{c}{4\alpha}) - \frac{\beta^{3}}{\alpha} + \Big(\phi
U\big(\frac{\gamma}{\alpha} - \frac{c}{4\alpha} +
\frac{\beta^{2}}{\alpha}\big) \Big) = 0,
\end{equation}
\begin{equation}
U\big(\frac{\gamma}{\alpha} - \frac{c}{4\alpha} \big)=
\kappa_{3}\frac{\beta^{2}}{\alpha},
\end{equation}
\
\\
Combining (2.21) with (2.17), (2.18) and (2.20) we obtain
\begin{equation}
\Big( \phi U\big(\frac{\gamma}{\alpha} - \frac{c}{4\alpha}\big)
\Big) = \frac{3\beta}{\alpha}\big[(\frac{\gamma}{\alpha} -
\frac{c}{4\alpha})^{2} - \frac{c}{4\alpha}\big].\end{equation}
 If $\gamma = \frac{c}{4}$ then the last relation yields
$\frac{3\beta c}{4} = 0$ which is a contradiction. Hence $\gamma
\neq \frac{c}{4}.$\qed

\begin{lemma}
Let M be a real hypersurface of a complex plane $M_{2}(c)$ $(c
\neq 0)$, satisfying $\phi l = l \phi$ on $ker(\eta)$. Then,
$\kappa_{3} = 0$ on $V'$.\end{lemma}
 Proof.\\
Because of (2.3), (2.4), (2.19), (2.22) and (2.23), the well known
relation $[U, \phi U] = \nabla_{U}\phi U - \nabla_{\phi U}U$ takes
the form
$$
[U, \phi U](\frac{\gamma}{\alpha} - \frac{c}{4\alpha})=
$$
$$-\frac{\kappa_{2}\kappa_{3}\beta^{2}}{\alpha}
  - \kappa_{3}\beta(\frac{\gamma}{\alpha} - \frac{c}{4\alpha} + \frac{\beta^{2}}{\alpha}) -
  \frac{3\beta\kappa_{3}}{\alpha}\big[(\frac{\gamma}{\alpha} -
\frac{c}{4\alpha})^{2} - \frac{c}{4}\big] - \kappa_{3}\beta\
(\frac{\gamma}{\alpha} - \frac{c}{4\alpha})$$

\
\\
On the other hand (2.18), (2.20), (2.22) and (2.23) yield:\\
$$[U, \phi U](\frac{\gamma}{\alpha} - \frac{c}{4\alpha}) = U\Big(
\phi U(\frac{\gamma}{\alpha} - \frac{c}{4\alpha})\Big) - \phi
U\Big(U (\frac{\gamma}{\alpha} - \frac{c}{4\alpha})\Big) = $$

$$ \frac{3(U\beta)}{\alpha}\big[(\frac{\gamma}{\alpha} -
\frac{c}{4\alpha})^{2} - \frac{c}{4}\big] - \frac{3\beta
(U\alpha)}{\alpha^{2}}\big[(\frac{\gamma}{\alpha} -
\frac{c}{4\alpha})^{2} - \frac{c}{4}\big] +
\frac{6\kappa_{3}\beta^{3}}{\alpha^{2}}(\frac{\gamma}{\alpha} -
\frac{c}{4\alpha}) - \frac{\beta^{2}}{\alpha}(\phi U(\kappa_{3}))$$
$$+ \frac{2\kappa_{3}\beta}{\alpha}(\frac{\gamma}{\alpha} -
\frac{c}{4\alpha})(\frac{\gamma}{\alpha} - \frac{c}{4\alpha} +
\frac{\beta^{2}}{\alpha}) - \frac{2\kappa_{3}\beta\gamma}{\alpha} -
\frac{\kappa_{1}\kappa_{3}\beta^{3}}{\alpha^{2}} -
\frac{\kappa_{3}\beta^{3}}{\alpha} -
\frac{3\kappa_{3}\beta^{3}\gamma}{\alpha^{3}}$$
$$+ \frac{3\kappa_{3}c \beta^{3}}{4\alpha^{3}}
$$
\ \\
The last equations using (2.15), (2.16) and (2.19) yield
\begin{equation}
\frac{3}{\alpha}[(\frac{\gamma}{\alpha} - \frac{c}{4\alpha})^{2} -
\frac{c}{4}](\xi\beta) -
\frac{3\beta}{\alpha^{2}}[(\frac{\gamma}{\alpha} -
\frac{c}{4\alpha})^{2} - \frac{c}{4}](\xi\alpha) -\beta(\phi U
\kappa_{3}) =
\end{equation}
$$
[2c - \beta \kappa_{2} + \frac{\beta^{2}}{\alpha} \kappa_{1} -
8(\frac{\gamma}{\alpha} - \frac{c}{4\alpha})^{2} -
\frac{5\beta^{2}}{\alpha}(\frac{\gamma}{\alpha} -
\frac{c}{4\alpha})]\kappa_{3}
$$

Following a similar way, from the action of $[\xi, \phi U]$ on
$\frac{\gamma}{\alpha} - \frac{c}{4\alpha}$ we obtain
\begin{equation}
\frac{3}{\alpha}[(\frac{\gamma}{\alpha} - \frac{c}{4\alpha})^{2} -
\frac{c}{4}](\xi\beta) -
\frac{3\beta}{\alpha^{2}}[(\frac{\gamma}{\alpha} -
\frac{c}{4\alpha})^{2} - \frac{c}{4}](\xi\alpha) -\beta(\phi U
\kappa_{3}) =
\end{equation}
$$[\gamma - (\frac{\gamma}{\alpha} - \frac{c}{4\alpha})^{2} - \frac{6\beta^{2}}{\alpha}
(\frac{\gamma}{\alpha} - \frac{c}{4\alpha})]\kappa_{3} $$

 Comparing (2.24) with
(2.25) and by making use of (2.17) we obtain
$$
\kappa_{3}\big[(\frac{\gamma}{\alpha} - \frac{c}{4\alpha})^{2} -
\frac{c}{4}\big] = 0
$$
Let us assume there is a point p on $V'$ such that $\kappa_{3} \neq
0$. Then because of the continuity of $\kappa_{3}$ there exists a
neighborhood W(p) around p such that $\kappa_{3} \neq 0$. This fact
and the last equation imply that $(\frac{\gamma}{\alpha} -
\frac{c}{4\alpha})^{2} = \frac{c}{4} $ on W(p). Differentiating the
last equation along $\xi$ and because of lemma 2.3 we obtain
$\xi(\frac{\gamma}{\alpha} - \frac{c}{4\alpha}) = 0$. Combining the
last equation with (2.19) we are led to $\kappa_{3} = 0$, which is a
contradiction. Therefore W(p) is empty and $\kappa_{3} = 0$ on
$V'$.\qed
\
\\
\section{Proof of theorem 0.3.}
We first prove the following:
\begin{lemma}
Let M be a real hypersurface of a complex plane $M_{2}(c)$ $(c \neq
0)$, satisfying $\phi l = l \phi$ and $\nabla_{\xi}l = \mu\xi$ on
ker$\eta$ or on span$\{\xi\}$. Then the structure vector field $\xi$
is principal on M.
\end{lemma}
Proof.\\
In what follows we work on $V'$. If $\nabla_{\xi}l = \mu\xi$ holds
on ker$\eta$, then using (2.1), (2.3), (2.4), (2.6) and (2.12) we
analyze $(\nabla_{\xi}l)\phi U = \mu\xi,$ $(\nabla_{\xi}l) U =
\mu\xi, $ and obtain $\mu + \beta\gamma = 0$ and $\mu = 0$. Since we
work on $V'$, the last two equations yield $\gamma = 0$ on $V'$. If
$\nabla_{\xi}l = \mu\xi$ holds on span$\{\xi\}$, then using (2.2),
(2.12) and (2.1) we also obtain $\gamma = 0$. From (2.20), (2.23)
and $\gamma = 0$ we obtain
\begin{equation}
\kappa_{1} = -4\alpha.
\end{equation}
The relations (2.18), (3.1) and $\gamma = 0$ yield
\begin{equation}
\phi U \beta = - \frac{c^{2}}{16\alpha^{2}} +
\frac{c\beta^{2}}{4\alpha^{2}} - 3\beta^{2}.
\end{equation}
The relations (2.23) and $\gamma = 0$ imply
\begin{equation}
\phi U \alpha = -3\alpha\beta + \frac{3\beta c}{4\alpha}.
\end{equation}

Now, from lemma 2.4, $\gamma = 0$ on $V'$, (2.22), (2.19), (2.15)
and (2.16) we have
\begin{equation}
(U\alpha) = (\xi\alpha) = (U\beta) = (\xi\beta) = 0
\end{equation}
which implies $$[U, \xi]\alpha = 0.$$
 On the other hand, by virtue of
(2.2), (2.3) and (3.1) we obtain
$$[U, \xi]\alpha = ( - \frac{c}{4\alpha} +
\frac{\beta^{2}}{\alpha} + 4\alpha)(\phi U \alpha).$$ The last two
equations imply
\begin{equation}(\frac{\beta^{2}}{\alpha} - \frac{c}{4\alpha} + 4\alpha)(\phi U
\alpha) = 0.
\end{equation}
We distinguish two cases: A)$\frac{\beta^{2}}{\alpha} -
\frac{c}{4\alpha} + 4\alpha \neq 0$ and B)$ \frac{\beta^{2}}{\alpha}
-
\frac{c}{4\alpha} + 4\alpha = 0$.\\
\ \\
A)Let $V_{1}$ be the set of points p of $V'$ such that $ -
\frac{c}{4\alpha} + \frac{\beta^{2}}{\alpha} + 4\alpha \neq 0$.
Because of the continuity there exists a neighborhood of p on which
$ - \frac{c}{4\alpha} + \frac{\beta^{2}}{\alpha} + 4\alpha \neq 0$
holds. This condition and (3.5) imply
$$\phi U \alpha = 0.$$
The last equation combined with (3.3) yields
\begin{equation}
\frac{c}{4} = \alpha^{2} > 0
\end{equation}
On the other hand calculating the Lie bracket $[U, \xi]\beta$ by
using (2.2), (2.3) and (3.4) we obtain $(\frac{\beta^{2}}{\alpha} -
\frac{c}{4\alpha} +
 + 4\alpha)(\phi U \beta) = 0$ which,
because of the condition $\frac{\beta^{2}}{\alpha} -
\frac{c}{4\alpha}  + 4\alpha \neq 0,$ leads to
$$\phi U \beta = 0.$$
From the last equation, (3.2) and (3.6) we obtain
\begin{equation}
  \frac{c}{4} = -2\beta^{2} =< 0
\end{equation}
which is a contradiction because of (3.6). Therefore the set
$V_{1}$ is empty.\\
\ \\
B)Let $V_{2}$ be the set of points p of $V'$ such that
$\frac{\beta^{2}}{\alpha} - \frac{c}{4\alpha} + 4\alpha = 0$.
Because of the continuity there exists a neighborhood of p on which
$ \frac{\beta^{2}}{\alpha} - \frac{c}{4\alpha} + 4\alpha = 0$ holds.
The last equation can be written as $\frac{c}{4} = 4\alpha^{2} +
\beta^{2} = 0.$ Differentiating this equation with respect to $\phi
U$ we obtain $\beta(\phi U \beta) + 4\alpha(\phi U\alpha) = 0.$ This
relation because of (3.2) and (3.3) takes the form $c\alpha^{2} =
0$, which is a contradiction. Therefore $V_{2}$ is empty. From the
cases A) and B) we obtain that $V'$ is empty.  This means that M
consists of the sets $V_{0}$ ad $V''$ where $\xi$ is principal. So
$\xi$ is principal everywhere on M.\qed \
\\
\
\\
\textbf{Proof of theorem 0.3}\\
From lemma 3.1 we have on M:
\begin{equation}
A\xi = \alpha\xi ,\hspace{50pt} \alpha = g(A\xi, \xi)
\end{equation}
We consider a $\phi - basis$ $\big\{ V, \phi V, \xi\big\}$. From
(2.5) and (3.8) we obtain
\begin{equation}
l V = \frac{c}{4}V + \alpha A V, \hspace{20pt} l\phi V =
\frac{c}{4}\phi V + \alpha A\phi V.
\end{equation}
The relations (3.9) and $\phi l = l\phi$ imply
\begin{equation}
  (A\phi - \phi A)V = 0.
\end{equation}
 On the other hand the action of $\phi$ on the second of (3.9) yields
$\phi (l\phi V) = \frac{c}{4}\phi^{2}V + \alpha \phi A\phi V$,
which, by virtue of (1.1), is written in the form
\begin{equation}
(\phi l)\phi V = -\frac{c}{4}V + \alpha (\phi A)\phi V,
\end{equation}
Moreover, from (1.1) and (3.9) we obtain:
$$(l\phi) \phi V = l \phi^{2}V = - l V =
-\frac{c}{4}V  - \alpha A V,$$ therefore
\begin{equation}
(l\phi) \phi V = -\frac{c}{4}V + \alpha (A\phi) \phi V
\end{equation}
Comparing (3.11) and (3.12), and using $\phi l = l \phi$ we have
\begin{equation}
  (A\phi - \phi A)\phi V = 0.
\end{equation}
 But from (1.1) and (3.8) we also have
\begin{equation}
  (A\phi - \phi A)\xi = 0.
\end{equation}
Therefore, (3.10), (3.13) and (3.14) imply that $A\phi = \phi A$.
This
result and the theorems (0.1) and (0.2) complete the proof of the theorem 0.3. \qed \\
\
\\
\section{Proof of theorem 0.4}\
\
\\
We define the sets:\\
$V_{0}$ = \{$p\in M/\alpha = 0$ in a neighborhood of M\},\\
$V'$ = \{$p\in M/\alpha \neq 0$, $\beta \neq 0$ in a neighborhood of M\},\\
$V''$ = \{$p\in M/\alpha \neq 0$, $\beta = 0$ in a neighborhood of M\}.\\
\
\\
Obviously $\xi$ is principal in $V''$. In addition we have proved in
section 2 that in every hypersurface satisfying $\phi l = l \phi$,
$\xi$ is
principal on $V_{0}$. In what follows we work on $V'$.\\

\begin{lemma}
Let M be a real hypersurface of a complex plane $M_{2}(c)$ $(c
\neq 0)$, satisfying $\phi l = l \phi$ and $lA = Al$ on $ker(\eta)$
or on span$\{\xi\}$. Then the structure vector field $\xi$ is
principal on M.
\end{lemma}
Proof.\\
If $lA = Al$ holds on $ker(\eta)$, we have $lAU = AlU$. The last
equation, (2.1), (2.6) and (1.6) yield $\gamma = 0.$ If $lA = Al$
holds on span$\{\xi\}$, then $lA\xi = Al\xi$. The last equation with
(1.6), (2.12) and (2.1) also yield $\gamma = 0$. Now the rest of the
proof is similar to this one of lemma 3.1. \qed. \
\
\\
\ \\Theorem 0.4 is proved by virtue of lemma 4.1 and following a
similar proof of theorem 0.3.
\bibliographystyle{amsplain}

\end{document}